\UseRawInputEncoding

\documentclass[12pt,oneside,a4paper]{article}

\usepackage{amsfonts}
\usepackage{amssymb}
\usepackage{amsmath}

\input xy
\xyoption{all}

\newtheorem{theorem}{Theorem}[section]

\newtheorem{proposition}[theorem]{Proposition}
\newtheorem{defi}[theorem]{Definition}

\newtheorem{remark}[theorem]{Remark}
\newtheorem{prbl}[theorem]{Problem}
\newtheorem{exmpl}[theorem]{Example}

\newenvironment{definition}
{\begin{defi} \rm}{ \end{defi}}
\newenvironment{example}
{\begin{exmpl} \rm}{ \end{exmpl}}

\newcommand{\done}{\rule{2mm}{2mm} \vskip \belowdisplayskip}
\newcommand{\qed}{\hfill\done}

\newenvironment{proof}
{\begin{trivlist} \item[] { \indent \textsl{Proof}.}} {\qed
\end{trivlist}}

\begin{document}
\title{\textbf{Lightly chaotic functional envelopes}}

\author{Annamaria Miranda\\
Dept. Mathematics,\\
University of Salerno,\\
Italy.\\
\date{}
\quad e-mail: {\tt amiranda@unisa.it}}

 \maketitle

\begin{abstract}

 In this paper we introduce some weak dynamical properties by using subbases for the phase space.
 Among them, the notion of \textit{light chaos} is the most significant.
Several examples, which clarify the relationships between this
kind of chaos and some classical notions, are given. Particular
attention is also devoted to the connections between the dynamical
properties of a system and the dynamical properties of the
associated functional envelope.
  We show, among other things, that a
continuous map $f: X\rightarrow X$, where $X$ is a metric space,
is chaotic (in the sense of Devaney) if and only if the associated
functional dynamical system is lightly chaotic.



\end{abstract}

\medskip

\section{Introduction}

When describing topological properties in terms of open sets,
 it suffices to restrict attention to a fixed
 base and, sometimes, to a fixed subbase. As an example
the \textit{Alexander subbase Theorem} claims: "\textsl{Let $X$ be
a topological space and  $\mathcal{S}$ a subbase for $X$. Then $X$
is compact if and only if every open cover by members of
$\mathcal{S}$ has a finite subcover.}"

\noindent However, the equivalence is not obvious in general.
Except for some rare cases, expressing a given topological
property by elements of a fixed subbase defines a weaker property
than the previous one.  In this paper we introduce some 'subbasic'
dynamical properties and a weak form of chaos, the \textit{light
chaos}, giving several examples in order to shed some light on the
relationships between this form of chaos and the notions of
transitivity and dense periodicity. Many authors investigated
about the chaotic behavior of a dynamical system (see for example
\cite{AAB},  \cite{CH}, \cite{Ru}). It should be rather natural to
think that some chaotic behavior of a dynamical system $(X,f)$
reflects, in some way, to its functional envelope $(S(X), F_f)$,
and viceversa. Recent researches deal with these kind of
relationships(see \cite{AKS}, \cite{C}). We establish a clear
connection between a discrete dynamical system $(X,f)$ and the
associated "functional" system $(S(X), F_f)$ regarding the chaotic
behavior. We will show that a system $(X,f)$ is chaotic (in the
sense of Devaney) if and only if $(S(X), F_ f)$ is
$\textit{lightly chaotic}$. Moreover we will give results showing
the non-chaoticity of $(S(X), F_f)$ for many dynamical system
$(X,f)$.

Our motivation to study the connection between a chaotic dynamical
system and the chaotic properties of its functional envelope also
comes from  interesting applications. To study the dynamical
behavior of the functional envelope is, in a certain sense,
equivalent to study the dynamical behavior of the solution of an
hyperbolic differential equation. Research about dynamics on a
space of functions has been motivated by the importance that the
topic plays in the semigroup theory, in the theory of functional
difference, in the dynamical systems theory and, moreover, and in
the study of one dimensional wave equations (see \cite{AKS},
\cite{C}, \cite{LH}, \cite{R}).

A discrete dynamical system is a pair $(X,f)$ where $X$ is a
topological space (called \textit{phase space}) and $f:
X\rightarrow X$ is a continuous map (called $\textit{transition
function}$). We may associate to $(X,f)$ a discrete dynamical
system $(S(X), F_f)$ whose phase space is the set $S(X)$ of all
continuous self-maps on $X$, endowed with a suitable topology, and
the transition function $F_f: S(X)\rightarrow S(X)$ is defined by
$F_f (g) = f\circ g$ for every $g\in S(X)$. In other words, the
phase space of the system $(S(X), F_f)$ is formed by the
transition functions of all discrete dynamical systems having $X$
as phase space, while $F_f$ is related, in a natural way, to the
transition function of the dynamical system $(X,f)$.

\noindent The first fundamental step in this topic is to define a
function space dynamical system by endowing the phase set with a
suitable topology assuring the continuity for the transition
function. The most quoted in literature are the compact-open
topology and the point-open topology. The second one is to
investigate about its dynamical properties, in connection with
those satisfied by the original system.

In \cite{AKS}, Auslander, Kolyada and Snoha  first called \textit{
functional envelope} of $(X, f)$ the induced system $(S(X), F_f)$,
since  it always contains an isometric copy of $(X, f)$(see \cite
{AKS}, Proposition 2.4), and denoted it by $(S(X),F_f)$, an usual
notation in the semigroup theory, since $S(X)$ is a topological
semigroup under the composition and the compact-open topology.
Their paper \textit{Functional envelope of a dynamical
system}(\cite{AKS}) deals with the connections between the
properties of a system and the properties of its functional
envelope, with special attention to orbit closures, $\omega$-limit
sets, (non)existence of dense orbits and topological entropy.
Since $(X, f)$ can be viewed as a subsystem of $(S(X), F_f)$, it
inherits many dynamical properties from $(S(X), F_f)$, but the
converse depends on topology on $S(X)$. They showed that, when $X$
is a compact metric space and S(X) is equipped with the
compact-open topology, some properties can be carried over from
$(X, f)$ to $(S(X); F_f)$, but there are other properties such as
topological transitivity, weakly mixing, mixing and chaoticity,
that cannot be easily transferred to the functional envelope. In
particular, they studied the following question which appears in
semigroup theory. Let $X$ be a compact metric space and S(X) be
equipped with the compact-open topology. Are there two elements
$f$ and $\varphi$ in $S(X)$ such that the set $O_{F_f} (\varphi)$
is dense in $S(X)$? In other words, for a given system $(X, f)$
with $X$ being a compact metric space, does the functional
envelope $(S(X), F_f)$ have a dense orbit? They showed that the
functional envelope of the full shift on $A^\mathbb{N}$, where $A$
is a compact metric space, contains dense orbits, but, in general,
the answer is negative (\cite{AKS}, Theorem 5.6 and Proposition
5.7]). So it is natural to ask if, for a compact metric space $X$
and a continuous map $f : X \rightarrow X$, there exists a
suitable, obviously coarser, topology $\mathcal{T}$ on $S(X)$ such
that the functional envelope $(S(X), F_f )$ of $(X, f)$ has at
least one dense orbit. In the paper \textit{Functional envelopes
relative to the point-open topology} (\cite{C}), the authors Chen
and Huang study this question. They consider a locally compact
separable topological space $X$ and the point-open topology
(denoted by $\mathcal{P}$), also known as \textit{pointwise
convergence topology}, on $S(X)$. In particular, for any
continuous map $f :[0,1]\rightarrow [0, 1]$, the functional
envelope $(S([0, 1]),\mathcal{P}, F_f)$ of $([0, 1], f)$ has no
dense orbits (\cite{C}, Theorem 2.1). Therefore it is not chaotic.
This leads them to restrict the point-open topology on a subset
$A$ of $X$, that is, to consider the topology $\mathcal{P}_A$
generated by the subbase $S_A = \{[x, U]$ : $x \in A$ and $U$ is
an open subset of $X\}$, where $[x, U] =\{\varphi \in S(X) :
\varphi(x)\in U\}$. They investigate the chaotic behavior of the
$(S_A(X), F_f)$, for a countable dense subset A of X, in relation
to that of the dynamical system $(X, f)$, where $X$ is a locally
compact separable metric space. Is a dynamical property, such as
transitivity, minimality, strong mixing of a system $(X, f)$,
absorbed by its functional envelope $(S_A(X), F_f)$? If not, what
is a coarser property satisfied by it?
 If $A$ is
countable and $X$ is a locally compact separable metric space they
show that:

\textit{If $(X, f)$ is weakly mixing and A is a countable dense
subset of X, then the functional envelope $(S_A(X), F_f)$ of
$(X,f)$ has at least one dense orbit.}

\medskip

\section{Preliminaries}

 In this section we recall some basic definitions and results involving topological dynamics and function space
 topologies useful in the sequel.
 We refer to \cite{A}, \cite{De}, \cite{GP} and \cite{W} for
definitions and results not explicitly  given. By a discrete
dynamical system we mean a pair $(X, f)$ where X is a (usually
compact metrizable) topological space and $f: X \rightarrow X$ is
a continuous map. Given a point $x\in X$,
%its trajectory is the sequence $ x, f (x),
%f^2(x), . . .$ and
its orbit is the set $O_f(x)=\{x, f (x), f^2(x), . . .\}$. For
every $n\in \mathbb{Z}^+$, we define the iterates $f^n : X
\rightarrow X$ as $f^0 = id_X $(the identity map on X) and
$f^{n+1}= f^n \circ f$.  By a subsystem of a dynamical system $(X,
f)$ we mean a system $(Y, g)$ where Y is a closed $f$-invariant
subset of $X$ (i.e. $f(Y)\subset Y$) and g is the restriction of
$f$ to $Y$.

A dynamical system $(X, f)$ (or the map $f$) is called
\textit{(topologically) transitive} if for every pair of nonempty
open sets $U$ and $V$ there is a positive integer $k$ such that
$f^k(U)\cap V\neq\emptyset$.

Topological transitivity and the existence of a dense orbit are
not equivalent for $(X, f)$, (see \cite{KS}). It is easy to check
that any $T_1$ point-transitive system without isolated points is
transitive. Conversely if $X$ is separable metrizable and of
second category then transitivity implies the existence of a dense
orbit (see \cite{S}).

A point $x$ is \textit{periodic} if $f^k(x)=x$ for some integer
$k\geq 1$. The least $k$ for which this happens is called
\textit{period} of $x$. A dynamical system $(X, f)$ (or the map
$f$) is called \textit{periodically dense} if the set of periodic
points of $f$ is dense in $X$.

Now, let $(X, d)$ be a metric space. A continuous map $f:(X,
d)\rightarrow (X, d)$  \textit{exhibits sensitive dependence on
initial conditions}, in brief is \textit{sensitive},
 if there is some $\delta>0$ such that,
for every point $x$ and every neighborhood $V$ of $x$ in $(X, d)$,
there exist a $y\in V$ and an integer $k\geq 1$ such that
$d(f^k(x), f^k(y))\geq \delta$.

Transitivity, sensitivity and periodic points' density are
dynamical ingredients used to introduce various kind of
chaoticity. Among them we consider Devaney chaoticity.

A continuous map $f: (X, d)\rightarrow (X, d)$ is chaotic (in the
sense of Devaney)  if it is transitive, periodically dense and
sensitive (see \cite{D}). It is worth noting that if $X$ is
infinite, a continuous map $f: (X, d)\rightarrow (X, d)$ is
chaotic (in the sense of Devaney) if and only if it is transitive
and periodically dense. So, it is surprising that the sensitivity,
a metric property which plays a central role in the definition of
several kind of chaos,
 is, in the definition of Devaney  chaos, a redundant property (see \cite{B}).
% (see \cite{D} for the definition and \cite{B}
% for the result).
Moreover, if $X$ is an interval of the real line, then $f$ is
chaotic (in the sense of Devaney) if and only if it is transitive.
%All spaces considered here are assumed to be infinite and we refer
%to Devaney chaos when not explicitly specified.

Several topologies can be defined on the set $Y^X=\{f:X
\rightarrow Y\}$, where $X, Y$ are topological spaces (see
\cite{AD}, \cite{AKO},\cite{KR},\cite{MN}, \cite{M}, \cite{O1}).
Among them we shall consider the compact-open topology, the
point-open topology and the uniform convergence topology. Let
$\mathcal{K}[X]$ be the family of all compact subsets of $X$. For
each $K\in \mathcal{K}[X]$ and $G\in \mathcal{T}$ let us set

$$[K,G]=\{g\in Y^X: g(K)\subset G\}.$$

\noindent The family $\mathcal{S}=\{[K, G]: K\in \mathcal{K}[X],
G\in \mathcal{T} \}$ is a subbase for a topology $\mathcal{T}_{k}$
on $Y^X$, the \textit{compact-open topology}. Moreover, given a
point $x\in X$ and an open set $G\in \mathcal{T}$, let

$$[\{x\},G]=\{g\in Y^X: g(x)\in G\}.$$
\noindent The sets $[\{x\},G]$
 form a
subbase for a coarser topology $\mathcal{T}_{p}$ on $Y^X$, the
\textit{point-open topology} (or \textit{topology of pointwise
convergence}). Note that this is just the product topology on
$Y^X$. Evidently, $\mathcal{T}_{p}=\mathcal{T}_{k}$ iff every
compact subset of $X$ is finite. In particular this happens when
$X$ is discrete. Moreover, on equicontinuous families of functions
of $Y^X$ the compact-open topology reduces to the point-open
topology (see \cite{W}).

 Now, let $(Y,d)$ be a metric space and let
$\mathcal{T}_d$ be the topology on $Y$ generated by $d$. We may
define on $B(X, Y)$, the set of all continuous bounded maps from
$X$ to $Y$, the \textit{uniform metric}:

$$\widehat{d} (f,g) = sup_{x\in X} d(f(x),g(x)).$$

\noindent The topology generated by $\widehat{d}$,
$\mathcal{T}_{\widehat{d}}$, is called the \textit{uniform
convergence topology}. Let $C(X,Y)$ is the set of continuous maps
from $X$ to $Y$. If $X$ is compact then $B(X,Y)=C(X,Y)$ and the
topology $\mathcal{T}_{\widehat{d}}$  coincides with the
compact-open topology. Moreover, this is independent from the
choice of $d$, i.e., every metric on $Y$ equivalent to $d$
generates the compact-open topology. If $(X, d)$ is a compact
metric space, then $(C(X,X)$, $\mathcal{T}_{\widehat{d}})$ is a
separable complete metric space.

Let $(X, \mathcal{T})$ be a topological space. We will denote
$C(X,X)$, the set of all continuous self-maps on $X$, by using the
simpler notation $S(X)$ coming from semigroup theory (as in
\cite{AKS} and in \cite{C}). It is straightforward to check that
if $f: (X,\mathcal{T})\rightarrow (X,\mathcal{T})$ is a continuous
map, then the map

$$F_{f}:(S(X),\mathcal{T}_{k})\rightarrow
(S(X),\mathcal{T}_{k})$$

\noindent defined by $F_f(g)= f\circ g$ for every $g\in S(X)$ is
continuous. In this way we associate to any discrete dynamical
system $(X,f)$ the dynamical system $(S(X), F_f)$. We will
concisely denote the topological space $(S(X),\mathcal{T}_{k})$ by
$S_{k}(X)$. When $X$ is a compact metric space then $(S_{k}(X),
F_f)$ (considered as a metric space with the uniform metric or
with the Hausdorff metric applied to the graphs of maps) is called
$\textit{functional envelope}$ of $(X, f)$ (see Definition 1.1. in
\cite{AKS}). Let us note that one could use this name even in a
more general setting.

\begin{definition}
Let $(X, f)$ be a dynamical system given by a topological
Hausdorff space $X$  and a continuous map $f : X \rightarrow X$.
If $F_f: (S(X),\mathcal{T}_{S(X)})\rightarrow
(S(X),\mathcal{T}_{S(X)})$ is continuous for some topology
$\mathcal{T}_{S(X)}$ on $S(X)$, then we call the dynamical system
$(S(X), F_f)$ \textit{functional envelope} of $(X, f)$.
\end{definition}

\begin{proposition}
 If   \hbox {  }$\mathcal{T}_{S(X)}$ is the compact-open topology or the point-
open topology then the map
$F_f:(S(X),\mathcal{T}_{S(X)})\rightarrow
(S(X),\mathcal{T}_{S(X)})$ is continuous and the dynamical system
$(S(X), F_f)$ contains a subsystem topologically coniugate to the
original system $(X, f)$.
\end{proposition}
\noindent In particular, if $\mathcal{T}_{S(X})= \mathcal{T}_{p}$
we will denote $(S(X),\mathcal{T}_{p}$) by $S_p(X)$.

\noindent  It is easy to prove that the previous result holds when
$\mathcal{T}_{S(X)}$ is any $\lambda-open$ topology. The
$\lambda-open$ topologies (or set-open topologies) topologies, a
generalization of the compact-open and of the pointwise
convergence topologies, were first introduced by Arens and
Dugundji (see \cite{AD}).

 \noindent Sometimes it is also convenient to endow $S(X)$ with coarser topologies induced by reduced
 subbases. For example,
if $(X, f)$ is a dynamical system given by a separable space X
with metric $d$ and a continuous map $f : X \rightarrow X$,  $A =
\{a_1, a_2,....\}$ is a countable dense subset of $X$, then
$(S_A(X), F_f)$, where $S_A(X)$ is $S(X)$ endowed with the
point-open topology on the set $A$, is a functional envelope of
$(X, f)$ (see Proposition 4 in \cite{C}).

\medskip
\section{Light chaoticity}

Compactness is equivalent to compactness with respect to a
subbase, as the \textit{Alexander subbase Theorem} states.
 However, subbases are not always sufficient to describe a topological
property and this allow to introduce a weaker property.

Indeed, let $\emph{P}$ be a topological property. Let $X$ be a
topological space and $\mathcal{S}$ a subbase for $X$. We say that
$X$ is $\textit{lightly}-P$ with respect to $\mathcal{S}$ if $X$
satisfies $\emph{P}$ relatively to all subbasic open sets in
$\mathcal{S}$. Evidently, any topological space $\emph{P}$ is
lightly -\emph{ P}, but the converse is not true. Obviously, when
the equivalence doesn't occur it makes sense to consider a weaker
form of the property.

Our attention is devoted to some dynamical properties. We
introduce in particular the \textit{light transitivity}, the
\textit{light periodically density} and the \textit{light
sensitivity}, and define a weak form of chaoticity, the
\textit{light chaoticity}.

 Let $(X, f)$
be a discrete dynamical system and $\mathcal{S}$ a subbase for the
topological space $(X,\tau)$.
\begin{definition}  $(X, f)$ (or the map $f: (X, \tau)\rightarrow (X, \tau$) is said to be \textit{lightly transitive} (with
respect to $\mathcal{S}$), briefly
$L_\mathcal{S}$-\textit{transitive}, if  for every $U, V \in
\mathcal{S}-\{\emptyset\}$ there exists some positive integer $k$
such that $f^k(U)\cap V\neq\emptyset$.
\end{definition}

\begin{definition}  $(X, f)$ (or the map $f: (X, \tau)\rightarrow (X, \tau$) is said to be \textit{lightly periodically dense} (with
respect to $\mathcal{S}$), briefly
$L_\mathcal{S}$-\textit{periodically dense}, if every $U\in
\mathcal{S}-\{\emptyset\}$ contains a
    periodic point of $f$.
\end{definition}

\begin{definition}  $(X, f)$ (or the map $f: (X, \tau)\rightarrow (X, \tau$) is said to be \textit{lightly chaotic} (with
respect to $\mathcal{S}$), briefly
$L_\mathcal{S}$-\textit{chaotic}, if:
\begin{description}
    \item[LC1] $(X, f)$ is $L_\mathcal{S}$-transitive;
    \item[LC2] $(X, f)$ is $L_\mathcal{S}$-periodically dense.
\end{description}
\end{definition}
In other words, a continuous map is $L_\mathcal{S}$-chaotic if it
satisfies transitivity and dense periodicity restricted to some
subbase $\mathcal{S}$ for $X$.
 Evidently, any $L_\mathcal{S}$-chaotic map is a
$L_\mathcal{S'}$-chaotic map for each $S'\subset S$. Moreover, if
the topology  $\tau(S')$ generated by $\mathcal{S'}$ is strictly
weaker than $\tau (S)$ and $f': (X, \tau (S')\rightarrow (X,
\tau(S')$ is continuous, then the dynamical system $(X,f')$ is
lightly chaotic too.

It is clear from the definition that:
\begin{description}
    \item[(i)] Every transitive periodically dense map is
    $L_\mathcal{S}$-chaotic with respect to any subbase for $X$.

    \item[(ii)] Every transitive interval map is chaotic, and, a
    fortiori,  lightly chaotic.
    \item[(iii)] Since a transitive map $f: S^1\rightarrow S^1$ is
    chaotic if and only if it has a periodic point, it follows
    that, for transitive self-maps on $S^1$, chaos and light chaos
    coincide.
\end{description}

The condition (i) is not necessary.
\begin{example}
\textit{A lightly chaotic periodically dense map which is not
transitive}.

\noindent The map $f: \mathbb{R}\rightarrow \mathbb{R}$ defined by
$f(x)=-x$ is a non transitive periodically dense map. Moreover $f$
is lightly chaotic with respect to the subbase
$\mathcal{S}=\{]-\infty, a[, ]b, +\infty[\}_{a, b \in
\mathbb{R}}$. The map $f$ sends any half line in its opposite and
two half lines, both right(or left), are always not disjoint. So
$f(]-\infty, a[)\cap ]b, +\infty[\neq 0 $ and $f(]b, +\infty[)\cap
]-\infty, a[\neq 0 $ for every $a, b \in \mathbb{R}$.

 \end{example}

\begin{example}

 \textit{A lightly chaotic transitive map which is not
periodically dense}.

\noindent Consider the dynamical system $(\{0,1\}^\mathbb{N},
\sigma)$, where $\{0,1\}^\mathbb{N}$ is the Cantor set and
$\sigma$ is the shift map. Let $S$ be the set of all eventually
constant sequences in $\{0,1\}^\mathbb{N}$ and $s^*\notin S$ a
transitive point. Since $\overline {O_\sigma(s^*)}=
\{0,1\}^\mathbb{N}$ then the set $X= S \cup O_\sigma(s^*)$ is
dense in $\{0,1\}^\mathbb{N}$. Moreover $\sigma(X)\subset X$. Let
$s\in X$. If $s=\sigma^n(s^*)$, then
$\sigma(s)=\sigma^{n+1}(s^*)\in X$.  If $s\in S$ (that is to say
$s=(s_0, s_1,..., s_n,..)$ and there exists a positive integer $k$
such that $s_n$ is constant $\forall n\geq k$), then $\sigma
(s_n)$ is eventually constant. Now, let $g=\sigma_{|X}:
X\rightarrow X$. Since $\sigma$ is transitive and $\overline X =
\{0,1\}^\mathbb{N}$, then $g$ is transitive (or simply $\overline
{O_\sigma(s^*)}= X$). $S$ is closed in $X$, the periodic points of
$g$ are the constant sequences, so $g$ is not periodically dense.
We claim that $g$ is lightly chaotic with respect to the canonical
subbase given by the sets of the form $X\cap\prod_n D_n$ where
$|D_n|=2$ $\forall n\in \omega-\{k\}$ and $|D_k|=1$, for some $k$.
Let $G=X\cap\prod_n D_n$ with $|D_n|=2$ $\forall n\neq k$ and
$|D_k|=1$. We may assume that $D_k=\{0\}$, so $\overline
0=(0,0,....)$ is a periodic point of $g$ such that $\overline 0\in
G$. Therefore $g$ is lightly chaotic.
\end{example}

\noindent Periodically density doesn't suffice to ensure light
chaoticity as the following example shows.
\begin{example}
\textit{A  map  which is not lightly periodically dense nor
lightly transitive}.

 Let $f: \mathbb{R}\rightarrow \mathbb{R}$ given by
$f(x)=|x|$.  $f$ is not lightly chaotic with respect to any
subbase for $\mathbb{R}$. Let $\mathcal{S}$ be a subbase for
$\mathbb{R}$ endowed with the usual topology. We may assume that
$\mathcal{S}$ consists of intervals (if $\mathcal{S}'$ is the
family of all components of members of $\mathcal{S}$, then $f$ is
lightly chaotic with respect to $\mathcal{S}$ iff it is lightly
chaotic with respect to $\mathcal{S}'$).
 Now let $V_1, ... , V_k, S_1, ... , S_n \in \mathcal{S}$ such
 that $V_1  \cap ...\cap V_k =]-\infty, -1[$ and $ S_1 \cap ... \cap
 S_n = ]1, +\infty[$. We may assume that $V_i , S_j \neq \mathbb{R} $ $\forall
 i,j$. So $V_i = ]-\infty, a_i[$,
$S_j=]b_j, +\infty[$ $\forall
 i,j$. Then $a_i=-1$, $b_j=1$ for some $i$ and for some $j$. Therefore $]-\infty, -1[$,
 $]1, +\infty[\in  \mathcal{S} $, and $f^k (]1, +\infty[)\cap ]-\infty, -1[
 =0$ $\forall k$. Thus $f$ is not lightly transitive.
Moreover the set of periodic points $\mathbb{R}^+_0$ is evidently
never dense.
\end{example}

    Now, let $(X, d)$ be a metrizable space.

    \begin{definition} A continuous map $f: (X,
d)\rightarrow (X, d)$  is said to be \textit{lightly
sensitive}(with respect to a subbase $\mathcal{S}$), briefly
$L_\mathcal{S}$-\textit{sensitive},
 if there is some $\delta>0$ such that,
for every point $x$ and every subbasic neighborhood $V$ of $x$ in
$(X, d)$, there exist a $y\in V$ and an integer $k\geq 1$ such
that $d(f^k(x), f^k(y))\geq \delta$.
\end{definition}

 \noindent It is clear that, since $X$ is metrizable,  any chaotic
map (in the sense of Devaney) is also $L_\mathcal{S}$-
    sensitive.
\begin{example}
\textit{A non lightly sensitive map}.

\noindent Let $f: [-1, 1]\rightarrow [-1, 1]$ given by
$f(x)=\frac{x}{|x|+1}$.   $f$ is not lightly sensitive with
respect to any subbase for $\mathbb{R}$ endowed with the usual
topology. Consider $\delta>0$, the point $x=0$ and a subbasic
neighborhood $V$ for $x$ such that $diam(f(V))<\delta$. Then for
every $y\in V$ and every integer $k\geq 1$ we have $d(f^k(x),
f^k(y))=d(0, \frac{y}{k|y|+1})=
\frac{|y|}{k|y|+1}\leq\frac{|y|}{|y|+1}=d(f(0), f(y))< \delta$

%OSS. Non è lightly transitive nè lightly periodically dense (l'unico punto periodico è lo zero)
\end{example}

Recall that a dynamical system is chaotic (in the sense of
Devaney) if and only if it is transitive, sensitive and
periodically dense. If $X$ is infinite, any transitive and
periodically dense map $f: (X, d)\rightarrow (X, d)$ is sensitive
(see \cite{B}).   But the corresponding assertion for light
properties does not remain true as the following example shows.
\begin{example}
 \textit{A lightly chaotic map which is not
lightly sensitive}.

\noindent Let $C^+ \subset\mathbb{R}^3$ be the cone whose base is
the circle $S^1$ in the coordinate plane $xy$ and whose vertex is
$(0,0,1)$ and $C^- \subset\mathbb{R}^3$ its symmetric with respect
to the coordinate plane $xy$. Consider the surface $X=C^+ \cup
C^-$, and the map $f$ defined by $f((e^{2\pi i\theta}),
t)=(|t-1|e^{2\pi i (\theta + \alpha) }, -t)$, where $\alpha =
\frac{p}{q}$ and $p$ and $q$ are co-prime. Note that this is a
glissorotation by an angle $\alpha $. The dynamical system $(X,f)$
is lightly chaotic but it is not lightly sensitive. All the open
half-spaces in $\mathbb{R}^3$ are a subbase. Let $\mathcal{S}$ be
the induced subbase on $X$. Evidently, any subbasic open set
contains a periodic point since it contains a vertex, a fixed
point,  or intersects the double cone base in a periodic points of
period $q$. Moreover, if $U$ and $V$ are two subbasic open sets
then there is a positive integer $k$ such that $f^k(U) \cap V$.
Indeed, if $U$ contains a vertex, then $f(U)\cap V$. If $U$ does
not contain any vertex, then there is an integer $k\leq q$ such
that $f^k(U)\cap V$. However $f$ is not lightly sensitive.
Consider any $\delta>0$ and let $x\in X$ be a point having
altitude $0$. Now, if $U$ is a subbasic neighborhood such that
$diam (U)<\delta$, then $d(f^k(x), f^k(y))=d(0, y) < \delta$ for
every $y\in U$ and $k\in \mathbb{N}$.

\end{example}

As already noted, if $X$ is an infinite metrizable space then
every transitive periodically dense map is
    $L_\mathcal{S}$-chaotic and $L_\mathcal{S}$-
    sensitive with respect to any subbase for $X$, but there are dynamical
    systems satisfying both light chaoticity and light
    sensitivity
which are not chaotic, as the following examples show.
\begin{example}
\textit{A lightly chaotic, lightly sensitive map which is neither
transitive nor periodically dense or sensitive}.

 \noindent The
"truncated tent map by  $\frac{1}{2}$" $T_ \frac{1}{2}:
I\rightarrow I $ is
 defined by

$T_ \frac{1}{2}(x)=\left \{%
\begin{array}{lll}
    2x, & if & 0\leq x < \frac{1}{4}\\
    \frac{1}{2}, &  if & \frac{1}{4}\leq x < \frac{3}{4}\\
     -2x+2, &  if & \frac{3}{4}\leq x \leq 1\\
\end{array}%
\right.$

\noindent The map $f: I\rightarrow I $ symmetric of  $T_
\frac{1}{2}: I\rightarrow I $ with respect to the line
$y-\frac{1}{2}=0$, defined
$f(x)=\left \{%
\begin{array}{lll}
    -2x+1, & if & 0\leq x < \frac{1}{4}\\
    \frac{1}{2}, &  if & \frac{1}{4}\leq x < \frac{3}{4}\\
     2x-1, &  if & \frac{3}{4}\leq x \leq 1\\
\end{array}%
\right.$ is not periodically dense:  $\mathcal{U}=] \frac{1}{4},
\frac{1}{2}[$ is an open set such that $\forall x\in \mathcal{U}$
$f^k(x)=\frac{1}{2}\neq x$ $\forall k$. Moreover, $f$ is not
transitive. It suffices to observe that $f$ is not onto. Evidently
$f$ is not sensitive. If $x=\frac{1}{2}$ and
$U=]\frac{1}{2}-\epsilon, \frac{1}{2}+\epsilon[$ for some
$\epsilon<\frac{1}{4}$, then $|f^k(x)-f^k(y)|=0$ for every $y\in
U$.

It is straightforward to check that $f$ is lightly chaotic with
respect to the subbase $\mathcal{S}=\{[0,a[, ]b,1]\}_{0<a\leq 1,
0\leq b <1}$. Moreover, $f$ is lightly sensitive. Let
$0<\delta<\frac{1}{4}$, $x\in X$ and $U$ a neighborhood of $X$. If
$U=[0,a[$ where $0<a\leq 1$ then there exists some positive
integer $k\geq 1$ such that $|f^k(x)-f^k(0)| = |f^k(x)-1|\geq
\delta$. If $U=]b,1]$ where $0\leq b <1$ then there exists some
positive integer $k\geq 1$ such that $|f^k(x)-f^k(1)|=
|f^k(x)-1|\geq \delta$.
 \end{example}

It is worth noting that, starting from the notion of topological
sensitivity (\cite {F}), it would be interesting to study
$L_\mathcal{S}$-\textit{topologically sensitivity}, too.

\section{Lightly chaotic functional envelopes}

 In the study of the relationships between the dynamical properties of a system and its functional envelope it is evident that the dynamical behaviour of the functional envelope
 is more complicated than that of the original system.
In particular, the
 functional envelope of a chaotic dynamical system in general
 fails to be chaotic. It suffice to think to the transitivity, even for usual spaces (see, for example,
Corollary 5.5 and Theorem 5.6 in \cite{AKS})).

Our first result fits in the realm of these investigations.  Here
are some starting remarks.

\begin{remark}

 For every $f: I \rightarrow I$ continuous, $F_f: S(I)\rightarrow
S(I)$ is not chaotic. However, if $f$ is chaotic then $F_f$ is
lightly chaotic.

\noindent Let $f: S^1 \rightarrow S^1$ be a transitive continuous
map coniugate to an irrational rotation. Since $f$ has no periodic
points, $F_f$ is not lightly chaotic with respect to any subbase
for $S^1$.

\end{remark}

\noindent It emerges a certain difficulty in obtaining that the
functional envelope is  chaotic, even if the original system is
chaotic. So, a natural question arises: is the chaoticity for a
given dynamical system equivalent to the light chaoticity of its
functional envelope?

The answer is given by the following.
\begin{theorem}
Let  $(X,d)$ be a metric space, and  $f: (X.d)\rightarrow (X,d)$ a
continuous mapping. Then the following are equivalent
\begin{enumerate}
        \item $f$ is chaotic
    \item $F_f: S_k(X) \rightarrow S_k(X) $ is lightly
chaotic with respect to the canonical subbase.
       \item $F_f: S_p(X) \rightarrow S_p (X)$ is lightly
%\item $X$ is chaotic
%        \item $S_k(X)$ is lightly
%chaotic with respect to the canonical subbase.
%        \item $S_p(X)$ is lightly
chaotic with respect to the canonical subbase.
    \end{enumerate}
\end{theorem}

\begin{proof}
Let $\mathcal{S}_{k}=\{[K, G]: K\in \mathcal{K}[X], G\in
\mathcal{T}(d) \}$ be the canonical subbase for $(S(X),
\mathcal{T}_{k})$. Let $A=[F, U], B=[C, V]\in
\mathcal{S}-\{\emptyset\}$. We will show that $F_f^k (A)\cap B
\neq \emptyset$ for some positive integer $k$. Since $U,V$ are
non-empty open subsets of $(X, d)$, and $f$ is transitive, there
is some positive integer $k$ such that $f^k (U)\cap V \neq
\emptyset$. Let $q\in U$ such that $f^k(q)=p\in V$. Now, let $g\in
S(X)$ be the map defined by $g(x)=q$ for every $x\in X$. Observe
that $g\in A$, so $F_f ^k(g)\in F_ f ^k(A)$. Moreover, $(F_f^k
(g))(x) = F_{f^k} (g))(x)=(f^k \circ g)(x)=f^k (g(x))=f^k (q)=p\in
V$ for every $x\in X$. Therefore ${F_f}^k (g)\in B$. Hence $F_f^k
(A)\cap B \neq \emptyset$. Now let us show that every member of
$\mathcal{S}-\{\emptyset\}$ contains a periodic point of $F_f$.
Let $A=[F, U] \in \mathcal{S}-\{\emptyset\}$. Since $U$ is a
non-empty open subset of $X$ and $f$ is periodically dense, there
is some $k>0$ and $x_0\in U$ such that $f^k (x_0)=x_0$. Let  $g: X
\rightarrow X$ be the map given by $g(x)=x_0$ for every $x\in X$.
Then $g\in A$ and $(F_f^k (g))(x)= (F_{f^k} (g))(x)=(f^k \circ
g)(x)=f^k (g(x))=f^k (x_0)=x_0\in V$ for every $x\in X$, so $F_f^k
(g)=g$ and $g$ is a periodic point of $F_f$ contained in $A$.
Therefore $F_f$ is lightly chaotic. Then $(S(X), F_f)$ is lightly
chaotic with respect to any subbase of the point-open topology
contained in $\mathcal{S}_{k}$.

Conversely, we have to show that $f$ is transitive and
periodically dense. Let us check  first the transitivity. Let $U$
and $V$ be a pair of non-empty open subsets of $X$, and let us
pick some
%$x_0\in U$ and $y_0\in V$. Now $A=[\{x_0, y_0\}, U],
%B=[\{x_0, y_0\}, V] \in \mathcal{S}-\{\emptyset\}$.
$x_0\in X$. Then, denoted by $\mathcal{S}_{p}$ the canonical
subbase for $(S(X), \mathcal{T}_{p})$,
 $A=[\{x_0\}, U], B=[\{x_0\}, V] \in
\mathcal{S}_{p}-\{\emptyset\}$. Since $F_f$ is lightly chaotic
with respect to the canonical subbase $\mathcal{S}_{p}$, there is
some $k>0$ such that $F_f^k (A)\cap B \neq \emptyset$. So there
are two maps $g$ and $h$ such that $g\in A$, $h\in B$ and $F_f^k
(g)=h$. Now $g(x_0)\in U$ and $h(x_0)\in V$, so $f^k (g(x_0))=F_
{f^k}(g)(x_0)=(F_ f^k(g))(x_0)=h(x_0) \in f^k (U)\cap V \neq
\emptyset$. Therefore, $f$ is transitive.

\noindent Now let us show that $f$ is periodically dense. Let $U$
be a non-empty open subset of $X$ and let $x_0 \in U$. Since $F_f$
is lightly chaotic and $A=[\{x_0\}, U]\in
\mathcal{S}_{p}-\{\emptyset\}$, there is some $g\in A$ and $k>0$
such that $F_f^k (g) = g$. Thus $f^k
(g(x_0))=(F_f^k(g))(x_0)=g(x_0)\in U$, and this means that
$g(x_0)$ is a periodic point of $f$ contained in $U$. The map $f$
is periodically dense, hence chaotic.

\end{proof}

Let us recall also that a continuous self-map $f$ on a compact
metric space $(X,d)$ is said to be \textit{structurally stable} if
there is some $\epsilon > 0$ such that every $g\in S(X)$ with
$\widehat{d}(f,g)<\epsilon$ is topologically coniugate to $f$.

We introduce an useful notion.
\begin{definition}
Let  $f: X\rightarrow X$ onto, $f$ is \textit{onto-stable} if
there exists an open set $U$ of $S_{k}(X)$ such that $f\in U$ and
every $g\in U$ is onto.
\end{definition}

\begin{remark}
Every transitive structurally stable map  $f: X\rightarrow X$,
where $X$ is a compact (metric) space, is "onto-stable".
\end{remark}

\begin{remark}

    Let $f: X \rightarrow X$ be a structurally stable map
    where $X$ is compact. Then $f$ is \textit{onto-stable}.

\end{remark}

 The next result gives some additional informations  about the
chaotic behavior of $(S(X),F_f)$, where $X$ is a compact
metrizable space and $F_f: S_k(X) \rightarrow S_k(X)$ is a
continuous map. Let $\emph{P}(f)$ be the set of all periodic
points for $f$.
\begin{theorem}
 The
following hold:
\begin{description}
    \item[i)] Let $X$ be a (metric) first countable continuum and let $f: X\rightarrow X$ be a continuous
    map. If $|\emph{P}(f)|<c$ then $F_ f$ is not periodically dense.
    \item[ii)] If $f: I\rightarrow I$, where $I=[0,1]$, is periodically dense, then
    $F_f$ is not transitive.
    \item[iii)] Let $X$ is be a compact (metric) space and let  $f: X\rightarrow X$ be a continuous
    map. If $f$ is transitive structurally stable "onto-stable" then $F_f$ is not
    transitive.
    \item[iv)] Let $X$ is be a compact (metric) space and let $f: X\rightarrow X$ be a continuous
    map. If $f$ is "onto-stable" and $\emph{P}(f)\neq X$ then $F_f$ is not periodically dense.

\end{description}

\end{theorem}
\begin{proof}

\begin{description}
    \item[i)] Let us suppose that $g$ is a periodic point of $F_
    f$. We claim that $g$ must be constant. Let us take a positive
    integer $k$ such that ${F_f}^k(g)=g$, that is   $F_f^k(g)(x)=g(x)$  $\forall x\in
    X$, then $g(x)$ is a periodic point of $f$ $\forall x\in X$. So
    $g(X)\subset \emph{P}(f)$. By hypothesis $|g(X)|<c$. Since $g(X)$
    is a (metric) first countable continuum space, it follows that
    $|g(X)|=1$, namely $g$ is constant.
    Now, let $a,b\in X, a\neq b$ and let $U$ and $V$ two open
    sets of $X$ such that $a\in U$, $b\in V$ and $U\cap
    V=\emptyset$. Then $G = [\{a\}, U]\cap[\{b\}, V]$ is a
    non-empty open subset of $S(X)$ (since it contains the identity
    map), which does not contain constant functions, so $F_
    f$ is not periodically dense.

    \item[ii)] Since $S(I)$ is a separable complete metric
    space without isolated points, it is enough to show that $F_f$ has no dense orbit (see \cite{KS}), i.e.,
    $\overline {O_{F_f}(g)}\neq S(I)$ for every $g\in S(I)$, where $O_{F_f}(g)= \{ {F_f}^n(g): n\in
    \mathbb{N}\}=\{g, f\circ g, f^2\circ g,...\}$. If $g$ is a
    constant map, then $O_{F_f}(g)$ consists of constant
    functions. Now the set $V=[[0, \frac{1}{2}],]0, \frac{1}{4}[]\cap [\{1\},]\frac{2}{3},\frac{3}{4}[]$
is a non-empty open subset of $S(I)$, equipped with the
compact-open topology, which does not contain constant maps, so
$V\cap O_{F_f}(g)=\emptyset$. If $g$ is not constant, then $g(I)$
has non-empty interior. Since $f$ is periodically dense, $g(I)$
contains a periodic point $p$ of $f$. Let $q\in I$ such that
$g(q)=p$ and set $U= [I, I-O_{f}(p)]$. $U$ is a non-empty open
subset of $S(I)$ such that $U\cap O_{F_f}(g)=\emptyset$. In fact
$(f^n\circ g)(q)= f^n(g(q))= f^n(p)\in O_{f}(p)$ for every $n\in
\mathbb{N}$, so $f^n\circ g\notin U$ $\forall n$. Therefore $f$ is
not transitive.
    \item[iii)] Any transitive structurally stable map $f: X\rightarrow X$, where is $X$ is be a compact (metric) space is onto-stable (by Definition 2.3).
    Let $\mathcal{C}_0$ be the subspace of $S(X)$
    consisting of all onto maps.
    Since $f\in\mathcal{C}_0$, then $\mathcal{C}_0\neq\emptyset$.
    Moreover, $\mathcal{C}_0$ is closed. Indeed, if $g\in S(X)-\mathcal{C}_0$, let $a\in X-g(X)$ and let $U = [X, X-\{a\}]$.
    $U$ is open in $S(X)$, $g\in U$ and
    $U\cap\mathcal{C}_0=\emptyset$, so $g\notin cl\mathcal{C}_0$.
    Observe that $\mathcal{C}_0\neq X$. Moreover $\mathcal{C}_0$ is
    $F_f$-invariant: $\forall g\in \mathcal{C}_0, F_ f (g)=
    f\circ g$ is onto, so $F_f (g)\in \mathcal{C}_0$. Since $f$ is structurally stable then $Int
    \mathcal{C}_0\neq\emptyset$. In fact, there exists $\epsilon >0$ such
    that $B_{\epsilon}(f)$ consists only of maps coniugate to $f$.
    These are transitive maps, hence onto, so $B_{\epsilon}(f)\subset \mathcal{C}_0$. Therefore $F_f$ is not transitive.

\item[iv)] Let $X$ be a compact metric space and $f: X\rightarrow
X$ an "onto-stable" map. Let us prove that if $\emph{P}(f)\neq X$
then $F_f$ is not periodically dense. Let $\varepsilon>0$ be such
that $g$ is onto $\forall g\in B_\varepsilon(f)$. If $g$ is a
periodic point of $F_f$, then $\exists k>0$ such that ${F_f}
^k(g)=g$, that is $f^k(g(x))=g(x)$ $\forall x\in X$. So $g(x)\in
\emph{P}(f)$ $\forall x\in X$, and this means that $g(X)\subset
\emph{P}(f)$. Therefore, by hypothesis, it follows that $g(X)\neq
X$, so $g\notin B_{\epsilon}(f)$. Therefore $B_{\epsilon}(f)$ does
not contain periodic points of $F_f$.
\end{description}
\end{proof}

As we noted there exists a lightly sensitive lightly chaotic map
which is neither transitive nor periodically dense. This continues
to hold true in hypothesis of sensitivity.
\begin{example}
 \textit{A sensitive lightly chaotic map  which is neither transitive nor
periodically dense}.

\noindent Let  $f: I\rightarrow I $ be any chaotic map, e.g., the
tent map. Applying the previous Theorems 4.2 and 4.6, it follows
that $F_f: S(I)\rightarrow S(I)$ is neither transitive nor
periodically dense, but it is lightly chaotic. It remains to show
that $F_f$ exhibits sensitive dependence to initial conditions. We
must prove that $\exists \delta>0$ such that $\forall g\in S(I)$
and for all open sets  $U_g$ containing $g$, $\exists h\in U_g$
and $n>0$ such that $\widehat{\rho}({F_f}^n(g), {F_f}^n
(h))>\delta$. So, let $U_g=\cap_{i=1}^m [K_i, V_i]$ and take
$x_0\in K_1, g(x_0)\in V_1$. Since the map $f$ is chaotic and
therefore sensitive, called  $\delta$ the sensitive constant, then
$\exists y_0\in V_1$ and $n>0$ such that $\rho(f^n(g(x_0)),
f^n(y_0))>\delta$. Let $h: I\rightarrow I $ be a continuous map
such that $h(K_1)=y_0$ and $h_{|K_j}=f_{|K_j}$ $\forall j\geq 2$.
Now, $h\in U_g$. Moreover $\widehat{\rho}(\overline f^n(g),
\overline f^n(h))= \sup_{x\in I}\rho( f^n(g(x)), f^n(y))\geq \rho(
f^n(g(x_0)),
 f^n(y_0))>\delta$.
\end{example}

\section{Concluding remarks}

Future investigations could have two perspectives.

The first concerning the introduction of other \textit{light}
dynamical properties and the study of their interdependencies with
classical dynamical properties, the second one concerning their
connections with dynamical properties of the functional envelope.
 It might be worth considering set-open topologies and uniform convergence
topologies (see for example \cite{AKO}, \cite{KR}, \cite{MN},
\cite{M}, \cite{O1}).
 Moreover, analogously, it might be
interesting to study the connections when the envelope is an
hyperspace.

\end{document}